# SPECTRA OF RANDOM LINEAR COMBINATIONS OF MATRICES DEFINED VIA REPRESENTATIONS AND COXETER GENERATORS OF THE SYMMETRIC GROUP[1]

By Steven N. Evans

*University of California, Berkeley*

We consider the asymptotic behavior as $n \to \infty$ of the spectra of random matrices of the form
$$\frac{1}{\sqrt{n-1}} \sum_{k=1}^{n-1} Z_{nk} \rho_n((k, k+1)),$$
where for each $n$ the random variables $Z_{nk}$ are i.i.d. standard Gaussian and the matrices $\rho_n((k, k+1))$ are obtained by applying an irreducible unitary representation $\rho_n$ of the symmetric group on $\{1, 2, \ldots, n\}$ to the transposition $(k, k+1)$ that interchanges $k$ and $k + 1$ [thus, $\rho_n((k, k+1))$ is both unitary and self-adjoint, with all eigenvalues either $+1$ or $-1$]. Irreducible representations of the symmetric group on $\{1, 2, \ldots, n\}$ are indexed by partitions $\lambda_n$ of $n$. A consequence of the results we establish is that if $\lambda_{n,1} \geq \lambda_{n,2} \geq \cdots \geq 0$ is the partition of $n$ corresponding to $\rho_n$, $\mu_{n,1} \geq \mu_{n,2} \geq \cdots \geq 0$ is the corresponding conjugate partition of $n$ (i.e., the Young diagram of $\mu_n$ is the transpose of the Young diagram of $\lambda_n$), $\lim_{n\to\infty} \frac{\lambda_{n,i}}{n} = p_i$ for each $i \geq 1$, and $\lim_{n\to\infty} \frac{\mu_{n,j}}{n} = q_j$ for each $j \geq 1$, then the spectral measure of the resulting random matrix converges in distribution to a random probability measure that is Gaussian with random mean $\theta Z$ and variance $1 - \theta^2$, where $\theta$ is the constant $\sum_i p_i^2 - \sum_j q_j^2$ and $Z$ is a standard Gaussian random variable.

**1. Introduction.** Beginning with [31], there has been a large body of work on various aspects of the distribution of the eigenvalues of large real symmetric or complex self-adjoint random matrices. In most cases, the stochastic structure has been the simplest possible, namely, the entries of the

Received August 2007.
[1]Supported in part by NSF Grant DMS-04-05778.

matrix are independent apart from the constraints imposed by symmetry or self-adjointness. Simultaneously, there has been significant investigation of the eigenvalues of large Haar distributed matrices from the classical groups. Some of this research is surveyed in [5, 23].

Recently, there has been interest in random matrices with more structure, as evidenced by [2] on the spectra of large random Hankel, Markov and Toeplitz matrices which was initiated by a list of open problems in [1] (see also [13]). In this paper we investigate ensembles of highly structured random matrices defined in terms of representations of the symmetric group.

Let $\mathfrak{S}_n$ denote the *symmetric group* of permutations of the set $\{1, 2, \ldots, n\}$. If we write elements of $\mathfrak{S}_n$ in cycle notation and, following the usual convention, omit cycles of length one, then $(k, k+1)$ denotes the permutation that interchanges $k$ and $k+1$, while leaving all other elements of $\{1, 2, \ldots, n\}$ fixed. The *transpositions* $(k, k+1)$, $1 \leq k \leq n-1$, generate $\mathfrak{S}_n$ and, writing $e$ for the identity element of $\mathfrak{S}_n$, they satisfy the *Coxeter relations*

$$(1.1) \qquad (k, k+1)^2 = e,$$

$$(1.2) \qquad (k, k+1)(\ell, \ell+1) = (\ell, \ell+1)(k, k+1), \qquad |k - \ell| \geq 2,$$

$$(1.3) \quad (k, k+1)(k+1, k+2)(k, k+1) = (k+1, k+2)(k, k+1)(k+1, k+2).$$

Suppose that $\rho$ is a (finite dimensional) *unitary representation* of $\mathfrak{S}_n$. That is, $\rho$ is a homomorphism from $\mathfrak{S}_n$ into the group of $d \times d$ unitary matrices for some $d$ [in other words, $\rho(\pi)$ is a unitary matrix for all permutations $\pi$ and $\rho(\pi'\pi'') = \rho(\pi')\rho(\pi'')$ for any two permutations $\pi'$ and $\pi''$]. Note from (1.1) that the transposition $(k, k+1)$ is its own inverse, so that

$$\rho((k, k+1)) = \rho((k, k+1)^{-1}) = \rho((k, k+1))^{-1} = \rho((k, k+1))^*$$

and, hence, each of the matrices $\rho((k, k+1))$ is *self-adjoint*. Consequently, any matrix of the form $\sum_{k=1}^{n-1} c_k \rho((k, k+1))$ for real coefficients $c_k$ will also be self-adjoint and have real eigenvalues. In this paper we investigate the asymptotic behavior as $n \to \infty$ of this spectrum when the corresponding representations of $\mathfrak{S}_n$ are chosen appropriately and the coefficients are chosen at random as realizations of independent Gaussian random variables with mean 0 and variance $\frac{1}{n-1}$.

In order to describe our results, we need some basic facts from the representation theory of the symmetric group. This material is available from many sources, such as [7, 8, 18, 21, 22, 26, 27, 30].

First of all, any unitary representation may be decomposed into a direct sum of *irreducible* representations. This means that the corresponding matrices will be in block diagonal form with a block for each irreducible component. Since the measure that puts unit mass at each eigenvalue (counting multiplicities) of the random matrix we are considering will be just a mixture of the analogously defined measures for each block, it suffices to consider the case where the representation $\rho$ is irreducible.

Second, since we are only interested in spectra, all that matters to us is the equivalence class of the representation $\rho$: two unitary representations $\rho'$ and $\rho''$ are said to be *equivalent* if there is a unitary matrix $U$ such that $\rho''(\pi) = U^*\rho'(\pi)U$ for all $\pi \in \mathfrak{S}_n$.

The equivalence classes of irreducible representations of $\mathfrak{S}_n$ are indexed by the *partitions* of $n$: a partition of $n$ is a sequence $\lambda = (\lambda_1, \lambda_2, \ldots, \lambda_k)$ of integers such that $\lambda_1 \geq \lambda_2 \geq \cdots \geq \lambda_k > 0$ and $\lambda_1 + \lambda_2 + \cdots + \lambda_k = n$ (we write $\lambda \vdash n$ and $|\lambda| := n$). We will describe concretely a representation corresponding to each such partition in Section 4, but all that matters for us at the moment is the *dimension* of the representations in each equivalence class [the dimension of an irreducible representation $\rho$ is the dimension of the matrices $\rho(\pi)$]. A partition $\lambda = (\lambda_1, \lambda_2, \ldots, \lambda_k) \vdash n$ may be pictured as a *Young diagram* consisting of $n$ boxes arranged in $k$ left-justified rows, with $\lambda_i$ boxes in row $i$ (where we number the rows from the top to the bottom).

We assign coordinates to the boxes in the diagram of $\lambda$ using the convention for entries of a matrix: box $(i,j)$ is the box in the $i$th row from the top and the $j$th column from the left, where $1 \leq i \leq k$ and $1 \leq j \leq \lambda_i$. If $(i,j)$ is the coordinate of a box in the diagram of $\lambda$, we write $(i,j) \in \lambda$.

The *hooklength* corresponding to box $(i,j) \in \lambda$ is

$$h_{i,j} := |\{(i,j') \in \lambda : j' \geq j\} \cup \{(i',j) \in \lambda : i' \geq i\}|.$$

The dimension of the equivalence class of irreducible representations indexed by $\lambda \vdash n$ is given by the celebrated *hooklength formula*

$$(1.4) \qquad f^\lambda := \frac{n!}{\prod_{(i,j) \in \lambda} h_{i,j}}$$

(see, e.g., Section 3.10 of [26] for a discussion of the history of this result). An alternative formula for $f^\lambda$ is

$$(1.5) \qquad f^\lambda = \frac{n! \prod_{i<j}(\lambda_i - \lambda_j - i + j)}{\prod_i (\lambda_i + k - i)!}$$

(see, e.g., Exercise 9 in Section 4.3 of [8] or Example 1 in Section I.1 of [22]).

Given a partition $\lambda \vdash N$, there is another partition $\lambda' \vdash N$ called the *conjugate* of $\lambda$ whose Young diagram is the transpose of the Young diagram of $\lambda$; that is, $(i,j) \in \lambda'$ if and only if $(j,i) \in \lambda$.

Our principal result is the following.

THEOREM 1.1. *For $n \geq 1$, let $\lambda_n$ be a partition of some positive integer $N_n$. Let $\rho_n$ be an irreducible unitary representation of $\mathfrak{S}_{N_n}$ corresponding to* $\lambda_n$. Write $\Xi_n$ for the random probability measure that puts mass $1/f^{\lambda_n}$ at each of the $f^{\lambda_n}$ eigenvalues (counting multiplicity) of the random matrix

$$\frac{1}{\sqrt{N_n - 1}} \sum_{k=1}^{N_n - 1} Z_{n,k} \rho_n((k, k+1)),$$

where $Z_{n,1}, Z_{n,2}, \ldots, Z_{n,N_n-1}$ are independent standard Gaussian random variables. Suppose that $N_n \to \infty$ as $n \to \infty$ and

$$\lim_{n \to \infty} \frac{\sum_{(i,j) \in \lambda_n} (j - i)}{\binom{N_n}{2}} = \lim_{n \to \infty} \frac{\sum_i \binom{\lambda_{n,i}}{2} - \sum_j \binom{\lambda'_{n,j}}{2}}{\binom{N_n}{2}} = \theta$$

exists. Then $\Xi_n$ converges in distribution (with respect to the topology of weak convergence of probability measures on $\mathbb{R}$) to a random probability measure $\Xi_\infty$ that is Gaussian with random mean $\theta Z$ and nonrandom variance $1 - \theta^2$, where $Z$ is a standard Gaussian random variable. In particular, the (nonrandom) expectation measure $\mathbb{E}[\Xi_\infty]$ is standard Gaussian.

REMARK 1.2. It is clear that the conditions of the theorem hold for a sequence $\lambda_n$ of partitions if and only if they hold for the sequence $\lambda'_n$ of their conjugates, with the corresponding constant being $\theta' = -\theta$. Since the random variable $-Z$ has the same distribution as $Z$, we see that the distributions of the two limiting random measures are equal.

REMARK 1.3. If we think of both $\lambda_n$ and the conjugate $\lambda'_n$ for each $n$ as infinite, nondecreasing sequences of nonnegative integers by appending zeros, then the conditions of the theorem will hold with

$$\theta = \sum_i p_i^2 - \sum_j q_j^2$$

if

$$\lim_{n \to \infty} \frac{\lambda_{n,i}}{N_n} = p_i$$

and

$$\lim_{n \to \infty} \frac{\lambda'_{n,j}}{N_n} = q_j$$

for all $i, j$.

REMARK 1.4. As we observe in the proof of the theorem, the quantity

$$\frac{\sum_{(i,j) \in \lambda_n} (j - i)}{\binom{N_n}{2}} = \frac{\sum_i \binom{\lambda_i}{2} - \sum_j \binom{\lambda'_j}{2}}{\binom{N_n}{2}}$$

appearing in the statement of the theorem is just

$$\frac{\chi_n(1^{N_n-2}2^1)}{\chi_n(e)},$$

where $\chi_n$ is the character of the representation $\rho_n$, the notation $(1^{N_n-2}2^1)$ denotes the equivalence class of elements of $\mathfrak{S}_{N_n}$ that have $N_n - 2$ one-cycles and a single two-cycle, and $e \in \mathfrak{S}_{N_n}$ is the identity permutation. Suppose that $\lambda_n$ is chosen at random from the set of partitions of $N_n = n$ according to *Plancherel measure*—the probability measure that assigns mass $\frac{(f^\mu)^2}{n!}$ to each $\mu \vdash n$. It follows from the orthogonality relations for irreducible characters that $\frac{\chi_n(1^{n-2}2^1)}{\chi_n(e)}$ has mean 0 and variance $1/\binom{n}{2}$. Hence, by the Borel–Cantelli lemma, the condition of the theorem holds almost surely with $\theta = 0$. We note from [19] that

$$\sqrt{\binom{n}{2}} \frac{\chi_n(1^{n-2}2^1)}{\chi_n(e)}$$

is asymptotically Gaussian with mean 0 and variance 1. There has been considerable further work on the asymptotics of random character ratios: see, for example, [3, 9, 10, 11, 14, 17, 28, 29].

We give the proof of Theorem 1.1 in Section 2. We derive some combinatorial consequences of the proof of Theorem 1.1 in Section 3 and describe a particular concrete realization of the random matrices that appear in the theorem in Section 4.

**2. Proof of Theorem 1.1.** We begin with an indication of how the proof will proceed. For simplicity, we drop the index $n$ from our notation until further notice, so we are working with an irreducible unitary representation $\rho$ of $\mathfrak{S}_N$ associated with the partition $\lambda$ of $N$, and $\Xi$ is the random probability measure that puts mass $1/f^\lambda$ at each of the $f^\lambda$ eigenvalues of the random self-adjoint matrix

$$M := \frac{1}{\sqrt{N-1}} \sum_{k=1}^{N-1} Z_k \rho((k, k+1)).$$

The proof will verify that the moment generating function of $\Xi$ converges to that of the claimed distribution as $N \to \infty$.

The first step will be to note that the $s$th moment of $\Xi$ is given by $\int x^s \Xi(dx) = \frac{1}{f^\lambda} \operatorname{trace}(M^s)$. If $\chi$ is the *character* associated with the representation $\rho$ [i.e., $\chi(\pi) = \operatorname{trace} \rho(\pi)$ for $\pi \in \mathfrak{S}_N$], then the trace of $M^s$ is a normalized multiple sum of terms consisting of the product of $s$ of the Gaussian weights $Z_k$ multiplied by the value of $\chi$ evaluated at the product of $s$ of the transpositions $(k, k+1)$.

Recall that the value of $\chi(\pi)$ for $\pi \in \mathfrak{S}_N$ only depends on the *conjugacy class* to which $\pi$ belongs, and two permutations belong to the same conjugacy class if and only if they have the same *cycle type*: the cycle type of $\pi \in \mathfrak{S}_N$ is the partition of $N$ that lists the lengths of the cycles of $\pi$ in decreasing order. Recall also that $\chi(e)$, the value of $\chi$ at the identity permutation $e \in \mathfrak{S}_N$, is the dimension $f^\rho$ of $\rho$. We will show that the only terms that make a significant asymptotic contribution are those for which the product of the $s$ transpositions has cycle-type consisting of $N - 2r$ one-cycles and $r$ two-cycles for $0 \leq r \leq s$. The corresponding products of Gaussian weights are made up from $r + (s-r)/2$ distinct $Z_k$, with $r$ weights appearing once and $(s-r)/2$ appearing twice.

We will then observe that, for a fixed value of $r$ (and hence a fixed value of the character $\chi$), the normalized sum of the products of Gaussians converges to a *multiple Wiener integral* $\int_{[0,1]^r} W(dt_1) \cdots W(dt_r)$, where $W$ is a *standard white noise*. We will also recall classical expressions for the values of a character at permutations consisting of just one and two two-cycles and use them to prove that the value of $\chi$ evaluated at permutations consisting of $r$ two-cycles is asymptotically equivalent to $\theta^r f^\rho = \theta^r \chi(e)$.

Last, we will combine the connection between Wiener integrals and *Hermite polynomials* with the generating function for the Hermite polynomials to conclude that the limiting moment generating function of $\Xi$ is that of a Gaussian with the prescribed random mean and nonrandom variance.

Before proceeding to the details of the proof, we list the facts we will need about Hermite polynomials that may be found in [15, 24]. The Hermite polynomials $H_n$ are defined by the generating function

$$\sum_{n=0}^\infty t^n H_n(x) := \exp\left(tx - \frac{x^2}{2}\right).$$

Equivalently,

$$H_n(x) = \frac{(-1)^n}{n!} \exp\left(\frac{x^2}{2}\right) \frac{d^n}{dx^n} \exp\left(-\frac{x^2}{2}\right).$$

The polynomial $H_n$ has degree $n$ with leading coefficient $\frac{1}{n!}$. If $Z$ is a standard Gaussian random variable, then

$$\mathbb{E}[H_m(Z) H_n(Z)] = \begin{cases} \frac{1}{n!}, & \text{if } m = n, \\ 0, & \text{otherwise.} \end{cases}$$

If $W$ is the usual white noise on $[0,1]$ and $Z = W([0,1])$ (so that $Z$ is standard Gaussian), then $n! H_n(Z)$ is the *multiple Wiener integral* $\int_{[0,1]^n} W(dt_1) \cdots W(dt_n)$ for $n \geq 1$.

We now give the details of the proof of the theorem:

$$\int x^s \Xi(dx) = \frac{1}{f^\lambda} \operatorname{trace}(M^s)$$

$$= \frac{1}{f^\lambda} \frac{1}{(N-1)^{s/2}}$$

(2.1)
$$\times \sum_{k_1=1}^{N-1} \cdots \sum_{k_s=1}^{N-1} Z_{k_1} \cdots Z_{k_s}$$

$$\times \operatorname{trace}(\rho((k_1, k_1+1)) \cdots \rho((k_s, k_s+1)))$$

$$= \frac{1}{(N-1)^{s/2}} \sum_{k_1=1}^{N-1} \cdots \sum_{k_s=1}^{N-1} Z_{k_1} \cdots Z_{k_s} \frac{\chi((k_1, k_1+1)) \cdots (k_s, k_s+1))}{\chi(e)}.$$

We will repeatedly use the fact that, since the matrix $\rho(\pi)$ is unitary, $|\chi(\pi)| \leq f^\lambda = \chi(e)$.

We wish to decompose the sum (2.1) according to the cycle type of the product $(k_1, k_1+1) \cdots (k_s, k_s+1)$ and show that the contribution of most cycle types is negligible as $N \to \infty$. As a device for doing this, given a vector $k = (k_1, \ldots, k_s) \in \{1, \ldots, N-1\}^s$, we build an edge-labeled graph $G_k$ with vertices $\{1, \ldots, s\}$ by putting an edge between the vertices $1 \leq p < q \leq s$ if $|k_p - k_q| \leq 1$ and labeling this edge with the integer $k_q - k_p$. There are $4^{\binom{s}{2}}$ graphs on $s$ vertices where each edge is labeled with an element of $\{-1, 0, 1\}$, so the graphs that can arise from the construction are contained, independently of $N$, in a fixed finite set. Moreover, if $k'$ and $k''$ are two vectors with $G_{k'} = G_{k''}$, then, by the Coxeter relations (1.1–1.3), the product $(k'_1, k'_1+1) \cdots (k'_s, k'_s+1)$ has the same cycle type as the product $(k''_1, k''_1+1) \cdots (k''_s, k''_s+1)$ and so $\chi((k'_1, k'_1+1) \cdots (k'_s, k'_s+1)) = \chi((k''_1, k''_1+1) \cdots (k''_s, k''_s+1))$.

Consider first of all a graph $G^*$ which has at least one connected component of size 3 or greater and all other connected components of size at least 2. Since $\#\{k \in \{1, \ldots, N-1\}^s : G_k = G^*\} = o((N-1)^{s/2})$ and $\mathbb{E}[(Z_{k_1} \cdots Z_{k_s})^2]$ is uniformly bounded, it is clear that

(2.2) $$\lim_{N \to \infty} \frac{1}{(N-1)^{s/2}} \sum_{\{k : G_k = G^*\}} Z_{k_1} \cdots Z_{k_s} = 0, \quad \text{in } L^2.$$

Consider next a graph $G^*$ which has at least one connected component of size 3 and $r \geq 1$ connected components of size 1, so that the number of connected components of size 2 or greater is strictly less that $\frac{s-r}{2}$. Note that if $G_{k'} = G_{k''} = G^*$, and the vertex $p$ is one of the $r$ components of size 1 of $G^*$, then it is certainly the case that $k''_q \neq k''_p$ for $q \neq p$, and so if

(2.3) $$\mathbb{E}[Z_{k'_1} \cdots Z_{k'_s} Z_{k''_1} \cdots Z_{k''_s}]$$

is nonzero, then $k_p'' \in \{k_1', \ldots, k_s'\}$. Consequently,

$$\mathbb{E}\left[\left(\sum_{\{k:G_k=G^*\}} Z_{k_1}\cdots Z_{k_s}\right)^2\right]$$
$$= \mathbb{E}\left[\sum_{\{k',k'':G_{k'}=G_{k''}=G^*\}} Z_{k_1'}\cdots Z_{k_s'} Z_{k_1''}\cdots Z_{k_s''}\right]$$
$$= \mathrm{o}((N-1)^{r+(s-r)/2+(s-r)/2}),$$

and (2.2) holds for such a choice of $G^*$.

We are therefore left with considering the contribution of graphs with connected components of size at most 2.

We begin the analysis of such graphs by showing that the contribution from those with at least one edge labeled with a $+1$ or a $-1$ is negligible as $N \to \infty$. Suppose, therefore, that $G^*$ is a graph with connected components of size at most 2, $u$ edges labeled $+1$, $v$ edges labeled $-1$, and $w$ edges labeled 0, with either $u > 0$ or $v > 0$, so there $u + v + w$ connected components consisting of two vertices and $s - 2(u+v+w)$ connected components consisting of a single vertex. Consider vectors $k'$ and $k''$ with $G_{k'} = G_{k''} = G^*$. Consider first one of the $u+v$ components of size 2 of $G^*$ with an edge labeled $\pm 1$ consisting of the vertices $p$ and $q$ with $1 \le p < q \le s$. By the construction of $G_{k''}$ we have $k_q'' - k_p'' = \pm 1$ and $k_t'' \notin \{k_p'', k_q''\}$ for $t \notin \{p,q\}$. Thus, if (2.3) is nonzero, then for some $1 \le h \ne \ell \le s$ we have $k_h' = k_p''$ and $k_\ell' = k_q''$, and, by the construction of $G_{k'}$, this means that one of the $u+v$ components of size 2 of $G_{k'}$ with an edge labeled $\pm 1$ consists of the vertices $h$ and $\ell$. Next consider one of the components of size 1 of $G^*$ consisting of the vertex $r$. By the construction of $G_{k''}$, we have that

(2.4) $$|k_t'' - k_r''| > 1, \quad \text{for } t \ne r.$$

The expectation (2.3) will be zero unless for some $1 \le m \le s$ we have $k_m' = k_r''$. Moreover, we must have that the vertex $m$ is one of the components of size 1 of $G^*$: if $m$ belonged to a component of size 2 with its edge labeled 0, then (2.3) would still be zero because $\mathbb{E}[Z_{k_r''}^3] = 0$, and if $m$ belonged to a component of size 2 with its edge labeled $\pm 1$ with the other vertex being $h$, then $|k_h' - k_m'| = 1$ and $k_h' \ne k_g'$ for $g \ne h$ by construction of $G_{k'}$, and so (2.3) would still be zero because (2.4) would give $k_h' \ne k_t''$ for $1 \le t \le s$. Putting these observations together, we have

$$\mathbb{E}\left[\sum_{\{k',k'':G_{k'}=G_{k''}=G^*\}} Z_{k_1'}\cdots Z_{k_s'} Z_{k_1''}\cdots Z_{k_s''}\right]$$
$$= \mathrm{O}((N-1)^{u+v+w+(s-2(u+v+w))+w})$$

$$= \mathrm{O}((N-1)^{s-u-v})$$
$$= \mathrm{o}((N-1)^s),$$

and (2.2) holds for this choice of $G^*$.

We are now left with graphs that consist of $r$ connected components of size 1 and $t$ connected components of size 2 with an edge labeled 0, where $0 \le r \le s$ and $r + 2t = s$. If two such graphs $G^*$ and $G^{**}$ share the same value of $r$, then

$$\sum_{\{k:\, G_k = G^*\}} Z_{k_1} \cdots Z_{k_s} = \sum_{\{k:\, G_k = G^{**}\}} Z_{k_1} \cdots Z_{k_s} = \sum Z_{a_1} \cdots Z_{a_r} Z_{b_1}^2 \cdots Z_{b_t}^2,$$

where the sum is over distinct $1 \le a_1, \ldots, a_r, b_1, \ldots, b_t \le N-1$, and the cycle type of the product $(k_1, k_1 + 1) \cdots (k_s, k_s + 1)$ is the same whether $G_k = G^*$ or $G_k = G^{**}$, namely, $N - 2r$ 1-cycles and $r$ 2-cycles. The number of connected graphs on $s$ vertices with $r$ connected components of size 1 and $t$ connected components of size 2 is

$$\binom{s}{r}(2t-1)(2t-3)\cdots 1 = \binom{s}{r} \mathbb{E}[Z^{s-r}],$$

since there are $\binom{s}{r}$ to choose the vertices that will comprise the $r$ connected components of size 1 and $(2t-1)(2t-3)\cdots 1$ ways to match the remaining $2t$ vertices into unordered pairs that will comprise the $t$ connected components of size 2.

Reintroducing the index $n$, we have shown that in order to understand the joint asymptotic behavior as $n \to \infty$ of $\int x^s \, \Xi_n(dx)$, $s \ge 0$, it suffices to understand the joint asymptotic behavior of

$$\frac{1}{(N_n - 1)^{s/2}} \sum_{r=0}^{s} \binom{s}{r} \mathbb{E}[Z^{s-r}] \frac{\chi_n(1^{N_n - 2r} 2^r)}{\chi_n(e)} \sum Z_{n,a_1} \cdots Z_{n,a_r} Z_{n,b_1}^2 \cdots Z_{n,b_t}^2,$$

where $\chi_n$ is the character of the representation $\rho_n$, $\chi_n(1^{N_n - 2r} 2^r)$ is the value of this character on any permutation consisting of $N_n - 2r$ 1-cycles and $r$ 2-cycles, and the interior sum is over $1 \le a_1, \ldots, a_r, b_1, \ldots, b_t \le N_n - 1$, where $r + 2t = s$ and distance between each pair of indices is at least 2 (the fact that $\mathbb{E}[Z^{s-r}] = 0$ when $s - r$ is odd allows us to write the outer sum over all $0 \le r \le s$, including those for which $s - r$ is odd and so no such $t$ exists). Our estimates above show that we get the same asymptotic behavior if we replace the condition that the distance between each pair of indices is at least 2 by the the condition that the indices are distinct.

We may take all of the random variables $Z_{n,k}$, $1 \le k \le N_n - 1$, $n \ge 1$, to be defined on the same probability space by putting

$$\frac{1}{\sqrt{N_n - 1}} Z_{n,k} := W\left(\left[\frac{k-1}{N_n - 1}, \frac{k}{N_n - 1}\right]\right),$$

where $W$ is a white noise on $[0,1]$. It then follows from an argument similar to that in the proof of Proposition 1.1.2 of [24] that

$$\frac{1}{(N_n-1)^{s/2}} \sum Z_{n,a_1} \cdots Z_{n,a_r} Z_{n,b_1}^2 \cdots Z_{n,b_t}^2$$

converges in probability to

$$\int_{[0,1]^r} W(dt_1) \cdots W(dt_r) = r! H_r(Z)$$

as $n \to \infty$, where $Z = W([0,1])$ is standard Gaussian.

An expression for

$$\frac{\chi_n(1^{N_n-2} 2^1)}{\chi_n(e)}$$

is given in [6] (see also [16]) that is equal to both

$$\frac{\sum_i \binom{\lambda_{n,i}}{2} - \sum_j \binom{\lambda'_{n,j}}{2}}{\binom{N_n}{2}}$$

and

$$\frac{\sum_{(i,j) \in \lambda_n} (j-i)}{\binom{N_n}{2}}$$

(see Exercise 7 in Section I.7, Exercise 3 in Section I.1, and equation (1.6) of [22]).

Given any sequence of integers, there will be a further subsequence along which

$$\frac{\chi_n(1^{N_n-2r} 2^r)}{\chi_n(e)}$$

has a limit for all $r$ [recall that $|\chi_n(\pi)| \leq \chi(e)$]. Denote this subsequential limit by $\kappa_r$. By assumption, $\kappa_1 = \theta$. It follows that for each $z \in \mathbb{R}$ there is a probability distribution $\xi_z$ with the sequence of moments

$$\sum_{r=0}^{s} s(s-1) \cdots (s-r+1) \mathbb{E}[Z^{s-r}] \kappa_r H_r(z), \qquad s \geq 0.$$

Recall that the leading term of $H_n(z)$ is $\frac{z^n}{n!}$ and so

$$\lim_{z \to +\infty} \int \left(\frac{x}{z}\right)^s \xi_z(dx) = \kappa_s, \qquad s \geq 0,$$

so that the sequence $\kappa_s$ is itself the sequence of moments of a probability distribution, say, $\bar{\xi}$.

An expression for

$$\frac{\chi_n(1^{N_n-4}2^2)}{\chi_n(e)}$$

is given in [16] that, from Table 3 of [4], is equivalent to

$$\frac{4(N_n-4)!}{N_n!}\left(\left(\sum_{(i,j)\in\lambda_n}(j-i)\right)^2 - 3\sum_{(i,j)\in\lambda_n}(j-i)^2 + 2\binom{N_n}{2}\right).$$

It follows that $\kappa_2 \leq \kappa_1^2$ and, hence, by the Cauchy–Schwarz inequality, the probability distribution $\bar{\xi}$ is the point mass at $\kappa_1 = \theta$. Therefore, $\kappa_s = \theta^s$. Since this is true for any subsequential limit, we must have that

$$\lim_{n\to\infty}\int x^s \Xi_n(dx) = \sum_{r=0}^{s} s(s-1)\cdots(s-r+1)\mathbb{E}[Z^{s-r}]\theta^r H_r(Z)$$

in probability for all $s \geq 0$.

It remains to note that

$$\sum_{s=0}^{\infty}\frac{t^s}{s!}\sum_{r=0}^{s} s(s-1)\cdots(s-r+1)\mathbb{E}[Z^{s-r}]\theta^r H_r(z)$$
$$= \sum_{s=0}^{\infty}t^s \sum_{r=0}^{s}\frac{\mathbb{E}[Z^{s-r}]}{(s-r)!}\theta^r H_r(z)$$
$$= \left(\sum_{s=0}^{\infty}\frac{t^s}{s!}\mathbb{E}[Z^s]\right)\left(\sum_{s=0}^{\infty}(\theta t)^s H_s(z)\right)$$
$$= \exp\left(\frac{t^2}{2}\right)\exp\left(\theta tz - \frac{(\theta t)^2}{2}\right)$$
$$= \exp\left(\theta zt + (1-\theta^2)\frac{t^2}{2}\right),$$

which we recognize as the moment generating function of a Gaussian probability measure with mean $\theta z$ and variance $1 - \theta^2$.

REMARK 2.1. In the proof of the theorem we used formulae for the character ratio $\frac{\chi_n(1^{N_n-2r}2^r)}{\chi_n(e)}$ for $r = 1$ and $r = 2$ due, respectively, to [6] and [16]. While the Murnaghan–Nakayama rule that we will discuss in Section 3 gives a recursive procedure for evaluating such a ratio, closed form expressions for arbitrary character values appear to have only been known in a few cases until the recent work of [20] building on [4].

**3. Two combinatorial identities.** In the proof of Theorem 1.1, we showed that if $\lim_{n\to\infty} \frac{\chi_n(1^{N_n-2}2^1)}{\chi_n(e)} = \theta$ existed, then $\lim_{n\to\infty} \frac{\chi_n(1^{N_n-2r}2^r)}{\chi_n(e)} = \theta^r$ for all $r \geq 0$. The Murnaghan–Nakayama formula, which we review below, gives a recursive expression for such character ratios and, taking limits, we will derive a pair of apparently nontrivial combinatorial equalities.

Given a partition $\lambda = (\lambda_1, \lambda_2, \ldots, \lambda_k) \vdash n$ and a partition $\mu = (\mu_1, \mu_2, \ldots, \mu_\ell) \vdash p$ for $p < n$, we write $\mu \subset \lambda$ if $\ell \leq k$ and $\mu_i \leq \lambda_i$ for $1 \leq i \leq \ell$. In other words, the Young diagram of $\nu$ is a subset of the Young diagram of $\lambda$. We write $\mu \prec \lambda$ if $p = n - 2$, $\mu \subset \lambda$, and the difference between the Young diagrams consists of a connected "domino" made up of two boxes, that is, if $\{(i,j) \in \lambda\} \setminus \{(i,j) \in \mu\}$ is of the form $\{(i^*, j^*), (i^*, j^* + 1)\}$ (a horizontal domino) or $\{(i^*, j^*), (i^* + 1, j^*)\}$ (a vertical domino). In the former case set $\zeta(\mu, \lambda) := +1$, and in the latter case set $\zeta(\mu, \lambda) := -1$. Given partitions $\lambda \vdash n$ and $\nu \vdash n - 2m$ for some positive integer $m$, we extend the definition of $\prec$ by declaring that $\mu \prec \lambda$ if there is a chain of partitions $\mu = \nu_0 \prec \nu_1 \prec \cdots \prec \nu_m = \lambda$. In this case we extend the definition of $\zeta$ by putting

$$\zeta(\mu, \lambda) := \sum_{\mu = \nu_0 \prec \nu_1 \prec \cdots \prec \nu_m = \lambda} \prod_{h=0}^{m-1} \zeta(\nu_h, \nu_{h+1}),$$

where the sum is over all sequences of partitions $\mu = \nu_0 \prec \nu_1 \prec \cdots \prec \nu_m = \lambda$.

It follows from the Murnaghan–Nakayama formula (see Section 4.10 of [26], formula 2.4.7 of [18] or Example 5 in Section I.7 of [22]) that

$$\frac{\chi_n(1^{N_n-2r}2^r)}{\chi_n(e)} = \sum_{\mu \prec \lambda_n,\, \mu \vdash N_n - 2r} \zeta(\mu, \lambda_n) \frac{f^\mu}{f^{\lambda_n}}.$$

Fix a sequence of partitions $\lambda_n = (\lambda_{n,1}, \ldots, \lambda_{n,K})$, each with $K$ parts, such that $\lim_{n\to\infty} |\lambda_n| = \infty$ and

$$\lambda_{n,i} - \lambda_{n,i+1} = 2\eta_i + 1, \qquad 1 \leq i \leq K - 1,$$

for nonnegative integers $\eta_i$. The conditions of Theorem 1.1 hold with $\theta = \frac{1}{K}$.

Drop the index $n$ for the moment and consider partitions $\lambda = (\lambda_1, \ldots, \lambda_K)$ with $\lambda_i - \lambda_{i+1} = 2\eta_i + 1$, $1 \leq i \leq K - 1$, and $\mu = (\mu_1, \ldots, \mu_K)$ with $\lambda_i - \mu_i = 2\delta_i$ for nonnegative integers $\delta_i$ such that $\sum_i \delta_i = r$, so that $\mu \vdash |\lambda| - 2r$. Suppose, moreover, that $2\delta_i \leq 2\delta_{i+1} + 2\eta_i + 1$, $1 \leq i \leq K - 1$, or, equivalently,

(3.1) $$\delta_i \leq \delta_{i+1} + \eta_i, \qquad 1 \leq i \leq K - 1.$$

Then $\mu \prec \lambda$ and all sequences $\mu = \nu_0 \prec \nu_1 \prec \cdots \prec \nu_r = \lambda$ are obtained by at each stage suitably adding a horizontal domino.

We enumerate the admissible chains as follows. Define a partition $\alpha = (\alpha_1, \ldots, \alpha_K)$ by

$$\alpha_i := \delta_K + \eta_{K-1} + \cdots + \eta_i, \qquad 1 \leq i \leq K - 1$$

and

$$\alpha_K := \delta_K$$

(this partition may have some parts that are 0, so we are extending our previous definition slightly, but we do have $\alpha_1 \geq \alpha_2 \geq \cdots \geq \alpha_K \geq 0$). Define another partition $\beta = (\beta_1, \ldots, \beta_K)$ by

$$\beta_i := \alpha_i - \delta_i, \qquad 1 \leq i \leq K-1,$$

and

$$\beta_K := 0$$

[this partition has some parts that are 0, but we do have $\beta_1 \geq \beta_2 \geq \cdots \beta_K \geq 0$ from (3.1)]. Note that the Young diagram of $\beta$ is a subset of the Young diagram of $\alpha$. There is only one way to tile the set difference of the Young diagrams of $\lambda$ and $\nu$ with horizontal dominoes, and each such domino corresponds to a unique box in the set difference of the Young diagrams of $\alpha$ and $\beta$. Moreover, the successive addition of horizontal dominoes implicit in a chain $\mu = \nu_0 \prec \nu_1 \prec \cdots \prec \nu_r = \lambda$ correspond to a labeling of the $r$ boxes in the set difference of the Young diagrams of $\alpha$ and $\beta$ with the integers $1, \ldots, r$ such that the labels in any row or column are increasing: the box labeled $k$ corresponds to the $k$th domino to be added to produce the Young diagram of $\mu_k$ from the Young diagram of $\mu_{k-1}$. From Corollary 7.16.3 of [30], the number of such labelings is

$$\zeta(\mu, \nu) = r! \det\left(\frac{1}{(\alpha_i - \beta_j - i + j)!}\right)_{1 \leq i,j \leq K},$$

with the convention that the factorial of a negative integer is zero (in the usual terminology, this quantity is the number of standard skew Young tableaux with shape $\alpha/\beta$).

Returning to the sequence of partitions $\lambda_n$, put $\bar{\lambda}_i := \lambda_{n,i} - \lambda_{n,K}$ for $1 \leq i \leq K$ (these quantities are independent of $n$ by assumption). It follows from (1.5) that

$$r! \sum \det\left(\frac{1}{(\delta_j + 1/2(\bar{\lambda}_i - \bar{\lambda}_j - i + j))!}\right)_{1 \leq i,j \leq K}$$
$$\times \frac{\prod_{i<j}(\bar{\lambda}_i - \bar{\lambda}_j - 2\delta_i + 2\delta_j - i + j)}{\prod_{i<j}(\bar{\lambda}_i - \bar{\lambda}_j - i + j)!}$$
$$= \kappa_r K^{2r}$$
$$= K^r,$$

where the sum is over all nonnegative integers $\delta_1, \ldots, \delta_K$ with

$$\delta_i \leq \delta_{i+1} + \tfrac{1}{2}(\bar{\lambda}_i - \bar{\lambda}_{i+1} - 1), \qquad 1 \leq i \leq K-1$$

and $\sum_i \delta_i = r$.

Suppose now that we specialize the preceding discussion to the case where $\eta_1 = \eta_2 = \cdots = \eta_K = 0$. The partitions $\mu$ with $\mu \vdash |\lambda| - 2r$ have $0 \leq \delta_1 \leq \cdots \leq \delta_K$ and the computation of $\zeta(\mu, \lambda)$ reduces to finding the number of standard Young tableaux with shape $(\delta_K, \delta_{K-1}, \ldots, \delta_1)$ (i.e., the number of ways of labeling the $r$ boxes in the Young diagram corresponding to this partition with the integers $1, \ldots, r$ in such a way that the labels in any row or column are increasing). This latter quantity is just $f^{(\delta_K, \delta_{K-1}, \ldots, \delta_1)}$ and is given by the analogue of (1.5). After a little algebra, we get that

$$\frac{r!}{\prod_i (i-1)!} \sum \frac{\prod_{i<j}(\delta_j - \delta_i + j - i)^2}{\prod_i (\delta_i + i - 1)!} = K^r,$$

where the sum is over $0 \leq \delta_1 \leq \cdots \leq \delta_K$ with $\sum_i \delta_i = r$.

As a check on this result, consider the special case when $K = 2$. Then the sum is over $\delta_1 = q$, and $\delta_2 = r - q$, where $0 \leq q \leq \lfloor \frac{r}{2} \rfloor$, and, by the symmetry of the binomial distribution and the usual formula for its variance,

$$r! \sum_{q=0}^{\lfloor r/2 \rfloor} \frac{(r - 2q + 1)^2}{q!(r - q + 1)!}$$

$$= \frac{1}{r+1} 2^{r+1} 2^2 \frac{1}{2} \sum_{q=0}^{r+1} \binom{r+1}{q} \left(\frac{1}{2}\right)^{r+1} \left(q - \frac{r+1}{2}\right)^2$$

$$= \frac{1}{r+1} 2^{r+1} 2^2 \frac{1}{2}(r+1) \frac{1}{2} \frac{1}{2}$$

$$= 2^r,$$

as claimed.

**4. An explicit matrix realization.** *Young's orthogonal representation* is a concrete element $\rho$ of the equivalence class of irreducible unitary representations associated with a given partition $\lambda \vdash N$ (see, e.g., [12, 18, 25]). The representing matrices are real (and hence orthogonal). Their rows and columns are indexed by the $f^\lambda$ standard Young tableaux for the Young diagram of $\lambda$, and $\rho(\pi)$ has a particularly simple expression when $\pi$ is one of the Coxeter generators $(k, k+1)$.

To be more explicit, suppose that $T$ is a standard Young tableau for $\lambda \vdash N$ and $1 \leq k \leq N - 1$. Write $(i', j')$ [resp., $(i'', j'')$] for the coordinates of the box in which $k$ (resp., $k+1$) appears. Set

$$d(T, k) := j'' - i'' - j' + i'.$$

Denote by $(k, k+1)T$ the tableau formed from $T$ by interchanging $k$ and $k+1$—this will only be a standard tableau if $k$ and $k+1$ do not appear in the same row or the same column of $T$. The $(T,T)$ entry of $\rho((k,k+1))$ is $d(T,k)^{-1}$. The $(T,U)$ entry of $\rho((k,k+1))$ is $\sqrt{1-d(T,k)^{-2}}$ if $U = (k,k+1)T$ and 0 otherwise [so that each row and column of the symmetric matrix $\rho((k,k+1))$ has at most one nonzero off-diagonal entry].

For example, suppose that $\lambda = (N-1,1)$. A standard Young tableau for $\lambda$ will be one of the tableau $T_2, T_2, \ldots, T_N$, where $T_\ell$ has $\ell$ in the single box in the second row of $\lambda$ and the numbers $1, 2, \ldots, \ell-1, \ell+1, \ldots, N$ in order in the $N-1$ boxes of the first row. Note that $d(T_\ell, k) = 1$ if $\ell \notin \{k, k+1\}$, whereas $d(T_k, k) = k$ and $d(T_{k+1}, k) = -k$. Observe also that $(k, k+1)T_k = T_{k+1}$ and $(k, k+1)T_{k+1} = T_k$ if $k \neq 1$, whereas $(k, k+1)T_\ell$ is not standard if $k = 1$ or $k \neq 1$ and $\ell \notin \{k, k+1\}$. Consequently, the matrix $\rho((1,2))$ has $(T_2, T_2)$ entry $-1$, all other diagonal entries $+1$, and all off-diagonal entries 0. Similarly, the matrix $\rho((k,k+1))$ for $k \neq 1$ has $(T_k, T_k)$ entry $+\frac{1}{k}$, $(T_{k+1}, T_{k+1})$ entry $-\frac{1}{k}$, remaining diagonal entries $+1$, $(T_k, T_{k+1})$ and $(T_{k+1}, T_k)$ entries $\sqrt{1-\frac{1}{k^2}}$, and all other entries 0.

Returning to a general partition $\lambda$, it is straightforward to describe the variances and covariances of the various entries in the centered Gaussian random matrix

$$\frac{1}{\sqrt{N-1}} \sum_{k=1}^{N-1} Z_k \rho((k,k+1)),$$

where the $Z_k$ are i.i.d. standard Gaussian, but the formulae do not appear to simplify to anything particularly pleasant.

**Acknowledgment.** This paper was written while the author was a Visiting Fellow at the Mathematical Sciences Institute of the Australian National University. The author thanks the members of the Institute for their hospitality.

DEPARTMENT OF STATISTICS #3860
UNIVERSITY OF CALIFORNIA, BERKELEY
367 EVANS HALL
BERKELEY, CALIFORNIA 94720-3860
USA
E-MAIL: evans@stat.berkeley.edu